\documentclass[11pt]{article}

\usepackage{amsmath}
\usepackage{amssymb}
\usepackage{amsfonts}
\usepackage{amsthm}
\usepackage{amstext}
\usepackage{amsbsy}
\usepackage{amsopn}
\usepackage{amscd}
\usepackage{latexsym}
\usepackage{lscape}

 \newtheorem{theorem}{Theorem}[section]
 \newtheorem{corollary}[theorem]{Corollary}
 \newtheorem{lemma}[theorem]{Lemma}
 \newtheorem{proposition}[theorem]{Proposition}
 \newtheorem{observation}{Observation}
 \newtheorem{remark}[theorem]{Remark}
 \theoremstyle{definition}
 \newtheorem{definition}[theorem]{Definition}
 \newtheorem{example}[theorem]{Example}
 \newtheorem{question}{Question}

 \numberwithin{equation}{section}

 \newcommand{\Ga}{\Gamma}
 \newcommand{\ga}{\gamma}
 \newcommand{\Om}{\Omega}
 
 \newcommand{\Si}{\Sigma}
 \newcommand{\si}{\sigma}
 \newcommand{\De}{\Delta}
 \newcommand{\de}{\delta}
 \newcommand{\al}{\alpha}
 \newcommand{\be}{\beta}
 \newcommand{\vp}{\varphi}

 \newcommand{\ov}{\overline}
 \newcommand{\mc}{\mathcal}
 \newcommand{\mr}{\mathrm}

 \newcommand{\la}{\langle}
 \newcommand{\ra}{\rangle}
 \newcommand{\tm}{\times}
 
 \newcommand{\h}{\hat}
 \newcommand{\semi}{\rtimes}
 \newcommand{\sm}{\setminus}

 \newcommand{\es}{\emptyset}

 \newcommand{\sse}{\subseteq}
 \newcommand{\nsse}{\nsubseteq}

 \newcommand{\Gb}{\overline{G}}
 \newcommand{\Ab}{\overline{A}}
 \newcommand{\Bb}{\overline{B}}
 
 \newcommand{\T}{\mathcal{T}}

 \newcommand{\move}{\mathrm{move}}
 
 \newcommand{\Aut}{\mathrm{Aut}}

 \newcommand{\Sym}{\mathrm{Sym}}
 \newcommand{\Alt}{\mathrm{Alt}}

 \newcommand{\core}{\mathrm{core}}
 \newcommand{\Soc}{\mathrm{Soc}}
 
 \newcommand{\Fun}{\mathrm{Fun}}

 \def\neq{\not=}

\title{\textbf{On Triple Factorisations of Finite Groups}}

\author{S. Hassan Alavi and Cheryl E. Praeger\\
\\
School of Mathematics and Statistics,\\ The University of Western Australia, \\
AUSTRALIA\\
\small
\texttt{alavi@maths.uwa.edu.au} \hspace{.3cm}  \texttt{praeger@maths.uwa.edu.au}}
\date{August 20, 2009 
}
\begin{document}
\maketitle

\begin{abstract}
This paper introduces and develops a general framework for studying triple factorisations  of the form $G=ABA$ of finite groups $G$, with $A$ and $B$ subgroups of $G$. We call such a factorisation nondegenerate if $G\neq AB$. Consideration of the action of $G$ by right multiplication on the right cosets of $B$ leads to a nontrivial upper bound for $|G|$ by applying results about subsets of restricted movement. For $A<C<G$ and $B<D<G$ the factorisation $G=CDC$ may be degenerate even if $G=ABA$ is  nondegenerate. Similarly forming quotients may lead to degenerate triple factorisations. A rationale is given for reducing the study of nondegenerate triple factorisations to those in which $G$ acts faithfully and primitively on the cosets of $A$. This involves study of a wreath product construction for triple factorisations.
\end{abstract}

\begin{quote}
\end{quote}

\newpage
\section{Introduction} \label{sec-intro}

\noindent In this paper we initiate a general theory of triple factorisations $\mathcal{T}=(G,A,B)$
of the form $G=ABA$ for finite groups $G$ and subgroups $A,B$. A special case of such factorisations was introduced by Daniel Gorenstein \cite{Gor59} in 1959, namely independent $ABA$-groups in which every element not in $A$ can be written uniquely as $abc$ with $a,c\in A$ and $b\in B$ (see \cite{Gor59}, also \cite{Gor62,GorHers59}). Triple factorisations also play a fundamental role as \emph{Bruhat decompositions} in the theory of Lie type groups (see for example, \cite{Carter89}), and more generally in the study of groups with a $(B, N )$-pair (see \cite{Bourbaki02}). With the extra condition that $AB\cap BA=A\cup B$, Higman and McLaughlin \cite{HM62} showed in their
famous paper that the associated coset geometry is a
flag-transitive $2$-design, and that $G$ acts primitively on the points of this design (that is, $A$ is a maximal subgroup of $G$). More general factorisations of the form $G=ABC$ arise, for example, as \emph{Iwasawa decompositions} of semisimple Lie groups (see \cite{FS01,Iwasawa49}). Here we focus on the case $A=C$.

Special cases of a triple factorisation $\T=(G,A,B)$ occur if
one of $A, B$ is equal to $G$, in which
case $\T$ is said to be \emph{trivial}, or more generally, if
$G$ factorises as $G=AB$, and here we say that $\T$ is
\emph{degenerate}. If $G\ne AB$ we call $\T$
\emph{nondegenerate}. The term `degenerate' may seem a misnomer, as group factorisations $G=AB$ arise in many problems involving symmetry in
algebra and combinatorics (see \cite{AFG92,LPS90} and the more than 100
references to these monographs recorded in MathSciNet). As the theory of
factorisations $G=AB$ is well developed, both for finite soluble groups
and for finite almost simple groups $G$, for the purposes of the present
theory we regard them as degenerate. Throughout the paper $G$ will
denote a finite group.

It is tempting to replace proper subgroups $A,B$ in a nondegenerate
triple factorisation $\T=(G,A,B)$ by `maximal overgroups', that
is, by subgroups $C, D$ where $A\leq C, B\leq D$, and $C, D$ are maximal
subgroups of $G$. This guarantees that $(G, C,D)$ is a triple
factorisation, but it may be degenerate. Indeed, in Section~\ref{sec-quo-lift}, we
show by a simple example that there may be no choices of maximal
overgroups $C,D$ that give a nondegenerate \emph{lift} $(G,C,D)$, even with $G$
a simple group, see Example~\ref{ex-lift}. One of the major questions we
address in this paper is:
\emph{under what conditions does such a lift $(G,C, D)$ of a
nondegenerate triple factorisation  $\T=(G,A,B)$ remain
nondegenerate?}

Important tools for studying triple factorisations come from permutation
group theory. For each proper subgroup $H$ of a group $G$, the group $G$
induces a transitive action by right multiplication on the set
$\Om_H=[G:H]$ of right cosets of $H$, and the kernel of this action is the
\emph{core} of $H$ in $G$, namely $\core_G(H)=\cap_{g\in G}H^g$.
Thus for a triple $\T=(G,A,B)$ with $A$ and $B$ proper
subgroups, there are two such transitive actions, on $\Om_A$ and
$\Om_B$.
In Section~\ref{sec-criteria} we give two necessary and sufficient
conditions for $\T$ to be a triple factorisation, one for each
of these actions.
The second criterion is in terms of the existence of a certain subset of
$\Om_B$ with `restricted movement' (see
Proposition~\ref{prop-second-crit}). Application of the classification result in \cite{NeuP96} yields a non-trivial
improvement to the trivially obtained upper bound
$|G|\leq|A|^2\cdot|B|/|A\cap B|$, and a generalisation to Theorem \ref{bound} is given in Theorem \ref{bound-gen}.

\begin{theorem}\label{bound}
If $\T=(G,A,B)$ is a nontrivial triple factorisation, then
\[ |G|\leq\frac{|A|^2|B|}{|A\cap B|^2}-\frac{|A|\cdot|B|}{|A\cap B|}+|B|\]
with equality if and only if $G$ is a flag-transitive collineation group
of a projective plane with point set $\Om_B$ and having the $A$-orbit
containing the coset $B$ as a line.
\end{theorem}

We remark that the finite flag-transitive projective planes, and also
their flag-transitive collineation groups are known explicitly, see
\cite{Kantor87}. It turns out that studying the $G$-action on $\Om_A$
is even more fruitful.

Among the possible simplifications, or `reductions', in studying the
class of triple factorisations $\T=(G,A,B)$ are the formation
of quotients
and restrictions of factorisations, which may yield triple
factorisations for groups of order less than $|G|$. For a normal
subgroup $N$ of $G$, the \emph{quotient} of $\T$ modulo $N$ is
the triple $\T/N=(G/N, AN/N, BN/N)$. This is always a triple
factorisation, and as in the case of lifts, it may be trivial (if $N$ is
transitive on $\Om_A$ or on $\Om_B$), degenerate or nondegenerate.
In Section~\ref{sec-quo-lift} we discuss various necessary and/or sufficient
conditions under which $\T/N$ is nondegenerate. In particular,
we prove the following elementary and useful fact (see Lemma \ref{lem-quo1}(c)).

\begin{observation}\label{modcore}
If $N\leq \core_G(A) \core_G(B)$, then $\T$ is nondegenerate if and only if $\T/N$ is nondegenerate.
\end{observation}

If $\T$ is nondegenerate, then it is often sufficient to study
$\T/N$ with $N= \core_G(A)\core_G(B)$. That is to say, we may assume in
our analysis that $G$ acts faithfully on both $\Om_A$ and $\Om_B$, or in other words, $A$ and $B$ are \emph{core-free} in $G$ (see Remark \ref{rem-quo1}(a)).
Thus problems about triple factorisations may be expressed in the
language of transitive permutation groups.
In particular if some lift $(G,C,D)$ were nondegenerate, with $C$ or $D$
maximal in $G$, then the $G$-action on $\Om_C$ or $\Om_D$
respectively, would be primitive, and we could reduce our analysis to
the context of primitive permutation groups. Unfortunately, as we noted
above, such a reduction is not always possible simply by considering a lift.

We do however have a rationale for focusing on primitive permutation groups, and to explain this
we  study the $G$-action on $\Om_A$ in more detail. If
$\T=(G,A,B)$ is a nondegenerate triple factorisation and $A$ is
not maximal in $G$, then there is a subgroup $H$ such that $A<H<G$. This gives rise to a lift $\T^{'}=(G,H,B)$ of $\T$ and corresponding quotient $\T_1=\T^{'}/N$ (which is also a lift of $\T/N$) where $N=\core_G(H)$ (see Proposition~\ref{prop-wr1}). In Section~\ref{sec-over-A} we show
that the \emph{restriction} $\T|_H=(H,A, B\cap H)$ is also a triple
factorisation, and  moreover, that at least one of $\T_1$ and
$\T|_H$ is nondegenerate (see Lemmas~\ref{lem-quo1} and
\ref{lem-res1}). By Observation \ref{modcore}, we may assume
that $G$ induces a faithful and imprimitive action on $\Om_A$
preserving a block system $\Si$ corresponding to the right
$H$-cosets. By the Embedding Theorem for imprimitive permutation groups
(see Theorem~\ref{thm:embed}), we may assume further
that $G$ is a subgroup of the wreath product $G_0\wr G_1$ acting on
$\De\times\Si$, where $G_1\cong G/N$ is the group induced by $G$ on
$\Si$ and $G_0\cong H/{\core_H(A)}$ is the group induced by $H$ on the block $\De$
of $\Si$ containing $A$. Moreover, we construct a \emph{wreath triple factorisation}
$\T_0\wr\T_1$ of $G_0\wr G_1$, where $\T_0$ is the quotient of
$\T|_H$ modulo $\core_H(A)$ (see Definition~\ref{defn-wreath}).
Our main result is the following imprimitive reduction theorem. Its
proof requires an extension of the Embedding Theorem, which we prove in
Theorem~\ref{thm:embed}.

\begin{theorem}\label{thm:main}
Let $\T=(G,A,B)$ be a nondegenerate triple factorisation with
$\core_G(A)=1$, and suppose that $A<H<G$, and $G\leq G_0\wr G_1$
with $G_0=H/{\core_H(A)}$, $G_1=G/{\core_G(H)}$ as above. Then there are triple factorisations
$\T_0$ for $G_0$, $\T_1$ for $G_1$ and $\T_0\wr \T_1$ for $G_0\wr G_1$  such that $\T_0$ is a quotient of $\T|_{H}$,
$\T_1$ is a quotient of a lift of $\T$, $\T_0\wr \T_1$ is nondegenerate, and either
\begin{enumerate}
\item[(a)] $\T_1$ is nondegenerate, or
\item[(b)] the triple factorisations $\T_0$ and
$\T_0\wr\T_1$ are both nondegenerate, and the
restriction of $\T_0\wr\T_1$ to $G$ is a nondegenerate
lift $(G,A,D)$ of $\T$.
\end{enumerate}
\end{theorem}

We prove a more detailed version of Theorem~\ref{thm:main} in
Section~\ref{sec-wreath}. If $B$ is maximal in $G$, then in Theorem~\ref{thm:main}(b) the
restriction of $\T_0\wr\T_1$ to $G$ is equal to
$\T$ (see Corollary~\ref{cor-wr4}).

\subsection{Rationale for primitive triple factorisations}

Observation~\ref{modcore} and Theorem~\ref{thm:main} provide a reduction pathway
to the study of nondegenerate triple factorisations
$\T=(G,A,B)$ with $A$ maximal in $G$ and both $A$ and $B$
core-free in $G$, as follows. We call such $\T$ \emph{primitive}. Corollary \ref{cor:main} below shows how each nondegenerate $\T$ is associated with at least one primitive triple factorisation.

\begin{corollary}\label{cor:main}
Let $\T=(G,A,B)$ be a nondegenerate triple factorisation with $\core_G(A)=1$. Then there exist subgroups $H$, $K$ such that $A\leq K< H\leq G$, $K$ is maximal in $H$, $G=HB$, and the quotient of $(H,K,BN\cap H)$ modulo $\core_H(K)$ is a primitive triple factorisation, where $N=\core_G(K)$. In particular, if $A$ is maximal in $G$ (so that $K=A$ and $H=G$), then $\T$ is primitive.
\end{corollary}

\emph{From this point of view the basic nondegenerate triple factorisations
$\T=(G,A,B)$  to study are the primitive ones, that is, those with $A$ and $B$ core-free in
$G$, and $A$ maximal in $G$.} For this study, we may therefore regard $G$ as a primitive
permutation group on $\Om_A$ with point stabiliser $A$, and study
intransitive subgroups $B$ with the property of Proposition~\ref{prop-first-crit}.

Often the study of primitive permutation groups $G$ focuses on the
action of its socle, and indeed, important (Lie type) triple factorisations arising
from the `Bruhat decomposition' often have $G$ simple. However
reduction of triple factorisations to proper normal subgroups is not straightforward. In our
final Section~\ref{sec-res-trans-normal} we discuss this problem, giving a rather technical sufficient condition for restriction.

\section{An extended Embedding Theorem} \label{sec-embed}

\noindent  In this section we give some preliminary definitions and results, and in particular prove an extension of the imprimitive embedding theorem appropriate for application to triple factorisations.

\subsection{Notation and basic definitions}

The set of all permutations of a set $\Om$ is the \emph{symmetric group} on $\Om$, denoted by $\Sym(\Om)$, and a subgroup of $\Sym(\Om)$ is called a \emph{permutation group} on $\Om$. Also we denote by $\Alt(\Om)$ the \emph{alternating group} on $\Om$. If a group $G$ acts on $\Om$ we denote the \emph{induced permutation group} of $G$ by $G^{\Om}$, a subgroup of $\Sym(\Om)$.  We say that $G$ is \emph{transitive} on $\Om$ if for all $\al$, $\be \in \Om$ there exists $g\in G$ such that $\al ^g=\be$. For a transitive group $G$ on a set $\Om$, a nonempty subset $\De$ of $\Om$ is
called a \textit{block} for $G$ if for each $g \in G$, either $\De^g=\De$, or $\De^g\cap\De=\es$; in this case the set $\Si=\{\De^g| \ g\in G\}$ is said to be a \textit{block system} for $G$. The group $G$ induces a transitive permutation group $G^{\Si}$ on $\Si$, and the set stabiliser $G_{\De}$ induces a transitive permutation group $G_{\De}^{\De}$ on $\De$. If the only blocks for $G$ are the singleton subsets or the whole of $\Om$ we say that $G$ is \textit{primitive}, and otherwise $G$ is \emph{imprimitive}.

\subsection{The Embedding Theorem}

Let $G\leq \Sym (\Om)$ be a transitive permutation group, and let

\begin{itemize}
  \item[(a)] $\Si=\{\De_1,\ldots, \De_\ell\}$ be a block system for $G$ in  $\Om$, set $\De=\De_1$ and $\al\in \De_1$;
  \item[(b)] $G_1=G^{\Si}\leq \Sym(\Si)$, and $G_0=G_{\De}^{\De}\leq \Sym(\De)$.
\end{itemize}

The group $G_1$ is determined by $G$, but the group $G_0$ may depend on the choice of block $\De$ in $\Si$. According to the \emph{Embedding Theorem}, this dependence is only up to permutational isomorphism: the Embedding Theorem gives a    \emph{ permutation embedding} $(\vp,\psi)$ from $(G,\Om)$ into $(G_0\wr G_1, \De\times \Si)$, that is to say, a   monomorphism $\vp:G\longrightarrow G_0\wr G_1$ and a bijection $\psi: \Om\longrightarrow \De\times \Si$ such that, for all $\de \in \De_i\in \Si$ and all $g\in G$,
\[(\psi(\de))^{\vp(g)}=\psi(\de^g).\]

\begin{itemize}
  \item[(c)] Note that if $\al \in \De$, $A=G_{\al}$, $A_0=A^{\De}$ and $A_1:=(G_1)_{\De_1}$, then the proof of the Embedding Theorem in \cite[Theorem 8.5]{NeumannBook} constructs a permutation embedding $(\vp,\psi)$ such that $\vp(A)\leq \h{A}$, where $\h{A}=(A_0\times G_0^{\ell-1}) \semi A_1 $.
\end{itemize}

When studying triple factorisations $G=ABA$ we need information about the subgroup $B$ as well as $A$. The group $B_1=B^{\Si}$ is determined by $B$ but  the group $B_0=B_{\De}^{\De}$ induced on $\De$ in general varies according to the choice of $\De$. However, if $B_1$ is transitive, then $B_0$ is unique up to permutation isomorphism and we have the following refinement of the Embedding Theorem.

\begin{theorem}[The Extended Embedding Theorem] \label{thm:embed}
Let  $G\leq \Sym(\Om)$ be transitive and let $\Si$, $\De$, $G_1$, $G_0$, $A$, $A_1$, $A_0$, and $\h{A}$ be as in (a)-(c) above. Suppose also that $B\leq G$, and let $B_1=B^{\Si}$, and $B_0=B_{\De}^{\De}$. Then there is a permutation embedding $(\vp,\psi)$ of $(G,\Om)$ into $(G_0\wr G_1, \De\times \Si)$ such that $\vp(A)=\h{A}\cap \vp(G)$, and if $B_1$ is transitive then $\vp$ may be chosen such that $\vp(B)\leq (B_0\wr B_1)$.
\end{theorem}
\begin{proof}
We express the permutation embedding $(\vp,\psi)$ explicitly, following the treatment in \cite[Theorem 8.5]{NeumannBook}. Here $G_0\wr G_1$ is identified as the semidirect product $\Fun(\Si,G_0) . G_1$, where $\Fun(\Si,G_0)$ consists of all functions $f: \Si \longrightarrow G_0$ and acts on $\De\times \Si$ via
\[(\de,\De_i)^f=(\de^{f(\De_i)},\De_i) \]
for all $\de \in \De$ and $\De_i\in \Si$, and $G_1$ acts naturally on $\De\times \Si$ by
\[(\de,\De_i)^\si=(\de,\De_i^\si),\]
and normalises $\Fun(\Si,G_0)$. If $B_1=B^{\Si}$ is not transitive, then all the assertions follow from \cite[Theorem 8.5]{NeumannBook}. Assume now that $B_1$ is transitive. Then, for each $i$, there is an element  $t_i\in B$ such that $\De_1^{t_i}=\De_i$, and the family $\{t_i\}_{i=1}^{\ell}$ is  a right transversal for $H:=G_{\De_1}$ in $G$. For $g\in G$ and $\De_i\in \Si$, define $i^g$ by $\De_i^{g}=\De_{i^g}$. As in \cite[Theorem 8.5]{NeumannBook}, the map $\vp: G\longrightarrow G_0\wr G_1$ is defined by
\begin{equation}\label{eq-embdding}
\vp(g)=f_gg^{\Si}, \hbox { where } f_g(\De_i)=(t_igt_{i^g}^{-1})^{\De}.
\end{equation}
Note that for each $i$, by the definition of $t_i$, $t_igt_{i^g}^{-1}\in G_{\De_1}=H$, so $f_g(\De_i)\in G_0$. Hence $f_g\in \Fun(\Si,G_0)$. This map is proved in \cite{NeumannBook} to be a monomorphism with $\vp(G)\leq G_0\wr G_1$, $\vp(A)=\h{A}\cap \vp(G)$, and $(\vp,\psi)$ is proved to be a permutation embedding, where $\psi:\Om\longrightarrow \De\times \Si$ is given by
\[\psi(\de)=(\de^{t_i^{-1}},\De_i), \hbox{ for } \de \in \De_i.\]
Since we have chosen all $t_i$ to lie in $B$, we have, for each $g\in B$, that $t_igt_{i^g}^{-1}\in B\cap H$, so $f_g(\De_i)\in B_0$. Thus $f_g\in \Fun(\Si,B_0)$ and hence $\vp(g) \in \Fun(\Si,B_0).B_1=B_0\wr B_1$. Therefore  $\vp(B)\leq B_0\wr B_1$.
\end{proof}

Note that if $B_1$ is not transitive, then $(B\cap H^{t_i})^{\De_i}$ may be very different from  $B_0=(B\cap H)^{\De_1}$ for some $i\not = 1$.

\begin{example}
Let $G=A_5 \wr S_4$, $A=A_4\times (A_5\wr S_3)$, $B=(A_5\wr S_2)\times (A_4\wr S_2)$, and $A<H=A_5\times (A_5\wr S_3)$. Then $\T=(G,A,B)$ is a nondegenerate triple factorisation, and relative to the block system $\Si=\{\De_1,...,\De_4\}$ corresponding to $H$, $(B\cap H^{t_i})^{\De_i}$ on blocks $\De_1$, $\De_2$ is $A_5$, and on $\De_3$, $\De_4$ is $A_4$.
\end{example}

\section{Triple factorisations: two criteria} \label{sec-criteria}

\noindent In what follows, $G$ is a group, $A$ and $B$ are subgroups of $G$, and we consider the triple $\T=(G,A,B)$. Recall that $\T$ is called a \textit{a triple factorisation} of $G$ if $G=ABA$. We use $\Om_A$ and $\Om_B$ to denote the set of right cosets of $A$ and $B$, respectively, that is, $\Om_A:=\{Ag| \ g \in G\}$ and $\Om_B:=\{Bg| \ g \in G\}$, and we note that $G$  acts naturally on both $\Om_A$ and $\Om_B$ by right multiplication. We call these actions \textit{right coset actions}. We give two different criteria for $\T=(G,A,B)$ to be a triple factorisation  based on the actions of $G$ on $\Om_A$ and $\Om_B$, respectively. The first criterion is connected to incidence geometries and we call it the \emph{geometric interpretation} of a triple factorisation. The second criterion relates to the notion of restricted movement of a set under a group action and we call this the \emph{restricted movement interpretation} of a triple factorisation. We were told of the first criterion  by Jan Saxl and it is discussed in~\cite{GJ}, while to our knowledge the second interpretation has not appeared in the literature before.

To explain the geometric interpretation, for a group $G$ and proper subgroups $A, B$, call the elements of $\Omega_A$ and $\Omega_B$ points and lines, respectively and define a point $Ax$ and a line $By$ to be incident if and only if $Ax\cap By\neq\emptyset$. An incident point-line pair $(Ax,By)$ is called a flag. It follows from this definition that $G$ preserves incidence and acts transitively on the flags of this geometry. Moreover, $\T=(G,A,B)$ is a triple factorisation if and only if any two points lie on at least one line (see Lemma 3 in \cite{HM62}). With extra conditions on a triple factorisation $\T$ making it  a \emph{Geometric $ABA$-group},  this incidence geometry becomes a $2$-design (see Proposition 1 in  \cite{HM62}).

\begin{proposition}[Geometric criterion]\label{prop-first-crit}
Let $A$ and $B$ be proper subgroups of a group $G$, and consider the right coset action of $G$ on $\Om_A:=\{Ag| \ g \in G\}$. Set $\al:=A \in \Om_A$. Then $\T=(G,A,B)$ is a triple factorization  if and only if  the $B$-orbit $\al^B$  intersects nontrivially each $G_\al$-orbit in $\Om_A$.
\end{proposition}

\begin{proof}
Set $\De:=\al^B$. Suppose that $\T$ is a triple factorisation for $G$, and consider a $G_\al$-orbit $\Ga:=\be^A$ in $\Om_A$. Since $G$ is transitive on  $\Om_A$, there exists $x \in G$ such that $\be =\al^x$, and since $G=ABA$,  $x=abc$, for some $a,c\in A$ and $b\in B$. Therefore, $\be^{c^{-1}}=\al^{xc^{-1}}=\al^{ab}=\al^{b}\in \De $, and also $\be^{c^{-1}} \in \be^A=\Ga$. Thus $\De\cap \Ga\not = \es$. Conversely, suppose that $\De$ intersects each $A$-orbit in $\Om_A$ nontrivially, and  let $x \in G$. Let $\Ga$ be the $A$-orbit containing $\al ^x$. By assumption, $\Ga \cap \De$ contains a point $\be$. Since $\be \in \De$, $\be =\al^b$ for some $b\in B$, and since $\be \in \Ga$, $\be=\al^{xa}$ for some $a\in A$. So $\al=\al^{xab^{-1}}$, and hence $c:=xab^{-1} \in G_\al=A$, that is, $x=cba^{-1} \in ABA$.
\end{proof}

The second criterion characterises triple factorisations of $G$ in terms of a subset of $\Om_B$ having restricted movement. For a finite subset $\Ga$ of $\Om$, the \emph{movement} of $\Ga$ under the action of a group $G$ on $\Om$ is defined by $\move(\Ga):=\max_{g\in G} |\Ga^g\sm \Ga|$.  If $\mr{move}(\Ga)<|\Ga|$, then $\Ga$ is said to have  \emph{restricted movement}. In other words, $\Ga^g\cap \Ga \not= \es$ for all $g\in G$.

\begin{proposition}[Restricted movement criterion] \label{prop-second-crit}
Let $A$ and $B$ be proper subgroups of a group $G$. Consider the right coset action of $G$ on $\Om_B$, and set $\be:=B \in \Om_B$. Then, $\T=(G,A,B)$ is a triple factorization  if and only if  the $A$-orbit $\be^A$  has restricted movement.
\end{proposition}
\begin{proof}
Let $\Ga:=\be^A$. Suppose that $\T$ is a triple factorisation and let $x\in G$.  Then there exist $a,c\in A$ and $b\in B$ with $x=abc$. Let $\ga:=\be^{bc}$. Since $G_{\be}=B$, we have $\be^b=\be$. Thus $\ga=\be^{bc}=\be^{c}\in \Ga$. Also, since $\Ga=\be^{A}$, we have $\ga=\be^{bc}\in (\be^{A})^{bc}=(\be^{A})^{abc}=\Ga^{x}$. Therefore, $\ga \in \Ga\cap \Ga^x$, that is,  $\Ga$ has  restricted movement. Conversely, suppose that $\Ga$ has restricted movement, and $x\in G$. Then $\Ga\cap \Ga^x\not= \es$, and hence there exist $a,c\in A$ with $\be^a=\be^{cx}$. Thus  $cxa^{-1}\in G_{\be}=B$,  so there exists $b\in B$ such that $x=c^{-1}ba$. Therefore, $G=ABA$.
\end{proof}

\begin{remark}\label{rem-crit1}
Let $\T=(G,A,B)$ be a triple factorisation, and consider $G$ acting on $\Om_B$. Let $\be=B \in \Om_B$, and let $k=|\be^A|$ and $m=\move(\be^A)$. If $m=0$, then $\be^A=\Om_B$, that is to say, $G=AB$, and vice versa. Therefore, $\T=(G,A,B)$ is nondegenerate  if and only if  $\be^A\not = \Om_B$  if and only if   $\move(\be^A)>0$. Hence, by Proposition \ref{prop-second-crit},  $\T$ is nondegenerate  if and only if  $1\leq \move(\be^A)\leq k-1$.
\end{remark}

Since $|AB|=|A||B|/|A\cap B|$, the equality $G=ABA$ implies that $|G|\leq |A|^2|B|/|A\cap B|$. A consequence of Proposition~\ref{prop-second-crit} is an improvement on this upper bound for $|G|$, namely the upper bound in Theorem~\ref{bound} in the introduction and its generalisation in Theorem \ref{bound-gen} below. A $2-(v,k,\lambda)$ \emph{design} is a set $\Om$ of $v$ \emph{points} together with a set of $k$-element subsets of $\Om$, called \emph{blocks}, such that each pair of points is contained in $\lambda$ blocks. Such a design is \emph{symmetric} if it has exactly $v$ blocks. (Note that this usage of the term `block' is different from the blocks of imprimitivity introduced in Section \ref{sec-embed}.) We apply the main result of \cite{Praeger97-Movement} to prove Theorem \ref{bound-gen} below, Theorem \ref{bound} follows immediately from Theorem \ref{bound-gen}. The special case of \cite[Theorem 1.1]{Praeger97-Movement} needed to prove Theorem \ref{bound} is proved in~\cite{NeuP96}.

\begin{theorem}\label{bound-gen}
Suppose $\T=(G,A,B)$ is a nondegenerate triple factorisation, $\be=B \in \Om_B$,  $k=|\be^A|$ and $m=\move(\be^A)$. Then
\begin{equation}\label{eq-bound-gen}
 |G|\leq \frac{|B|}{|A\cap B|}\left( \frac{|A|^2-m|A\cap B|^2}{|A|-m|A\cap B|}\right)
\end{equation}
with equality  if and only if  the $G$-translates of $\be^A$ form the blocks of a symmetric $2-(\frac{k^2-m}{k-m}, k,k-m)$ design with point set $\Om_B$ admitting $G$ as a flag-transitive group of automorphisms. Moreover, Theorem~\ref{bound} holds.
\end{theorem}
\begin{proof}
Consider the action of $G$ on $\Om_B$, and let $k=|\be^A|$. Then $k=|A|/|A\cap B|$, and by \cite[Theorem 1.1]{Praeger97-Movement}, we have $|\Om_B|\leq (k^2-m)/(k-m)$ with equality  if and only if  the $G$-translates of $\be^A$ form the blocks of a symmetric $2-(\frac{k^2-m}{k-m}, k,k-m)$ design with point set $\Om_B$ admitting $G$ as a flag-transitive group of automorphisms. Since $|\Om_B|=|G|/|B|$, the inequality (\ref{eq-bound-gen}) follows immediately. To prove Theorem \ref{bound}, note that the function, $f(m)=(k^2-m)/(k-m)$ increases as $m$ increases. Since $\T$ is nondegenerate, it follows from Remark \ref{rem-crit1} that $1\leq m\leq k-1$, and so by choosing $m=k-1$, the inequality in Theorem \ref{bound} follows from $|\Om_B|\leq k^2-k+1$  with equality  if and only if  $G$ is a flag-transitive group of collineations of a projective plane with point set $\Om_B$ such that $\Ga$ is a line.
\end{proof}
\section{Isomorphisms of triple factorisations}\label{sec-iso}

\begin{definition} \label{defn-iso1}
Let $\T:=(G,A,B)$ and $\T^{'}:=(G',A',B')$ be triple factorisations. We say that $\T$ is \emph{isomorphic} to $\T^{'}$ and write $\T\cong \T^{'}$ if there exists an isomorphism $\vp: G\longrightarrow G'$ such that $\vp(A)=A'$ and $\vp(B)=B'$.
\end{definition}

Suppose that $\T=(G,A,B)$ is a triple factorisation. Then for each $\si \in \Aut(G)$, $G=A^{\si}B^{\si}A^{\si}$. However, in general  the factorisation property is not preserved if we apply different automorphisms to $A$ and $B$ (see Lemma \ref{lem-res2}(a)). Moreover, if $\T^{'}=(G',A', B')$ is another triple factorisation then $\T$ and $\T^{'}$ may not be  isomorphic even if $G\cong G'$, $A\cong A'$ and $B\cong B'$ (see Lemma \ref{lem-res2}(b)).

\begin{lemma}\label{lem-res2}
\begin{enumerate}
  \item [(a)] There exist infinitely many nondegenerate triple factorisations $\T=(G,A,B)$ such that for some $x\in G$, $(G,A,B^x)$ is not a triple factorisation.
  \item [(b)] There exist infinitely many nondegenerate triple factorisations $\T=(G,A,B)$ and $\T^{'}=(G,A,B^x)$ where $x\in G$, such that $\T$ and $\T^{'}$ are not isomorphic.
\end{enumerate}
\end{lemma}

In order to prove Lemma \ref{lem-res2}(a) we prove a general result about triple factorisations $(G,A,B)$ with $G$ acting $2$-transitively on $\Om_A$, namely Proposition \ref{prop-res}, and then in Example \ref{ex-part(a)} we give an explicit family of examples.

\begin{proposition} \label{prop-res}
Suppose that $G$ is a $2$-transitive permutation group on a finite set $\Om$ of size $n$, let $\al \in \Om$, $A=G_{\al}$, and  $B<G$. Then  $\T=(G,A,B)$ is a triple factorisation  if and only if  $B \not \leq A$, and in this case $\T$ is nondegenerate  if and only if   $B$ is intransitive on $\Om$.
\end{proposition}
\begin{proof}
Since $G$ is $2$-transitive on $\Om$,  the $A$-orbits on $\Om$ are $\{\al\}$ and $\Om\sm \{\al\}$. Therefore, by Proposition \ref{prop-first-crit}, for a subgroup $B$ of $G$,  $\mc{T}=(G,A,B)$ is a triple factorisation   if and only if  $\al^B\not=\al$. Also $\mc{T}$ is nondegenerate  if and only if  $B$ is intransitive.
\end{proof}

\begin{example}\label{ex-part(a)}
Let $G=\Sym(\Om)$ with $\Om:=\{1,\ldots, n\}$, $n\geq 3$, $A:=G_1$, $B:=\la(1,3)\ra$, and $x:=(1,2)$. Then by Proposition \ref{prop-res}, $\T=(G,A,B)$ is a triple factorisation but $\T^{'}=(G,A,B^x)$ is not.
\end{example}

Similarly Lemma \ref{lem-res2}(b) is established by the next family of examples.

\begin{example}\label{ex-part(b)}
Let $G=\Sym(\Om)$ with $\Om:=\{1,\ldots, n\}$, $n\geq 5$. Let $A:=G_n$, $B:=\la(1,2), (3,4,n)\ra$, and $x:=(1,2,n)$. Then both $\T=(G,A,B)$ and $\T^{'}=(G,A,B^x)$ are  nondegenerate triple factorisations by Proposition \ref{prop-res}. However, $\T\not \cong \T^{'}$ for if $\vp \in \Aut(\Sym(\Om))$ with $\vp(A)=A$, then $\vp \in \Sym(\Om) \cap A=A$, and then $\vp(B)=\la (1^{\vp},2^{\vp}), (3^{\vp},4^{\vp},n)\ra$ whereas $B^{x}=\la (2,n), (3,4,1)\ra$. Since the $B^x$-orbit containing $n$ is $\{2,n\}$ while the $\vp(B)$-orbit containing $n$ contains $\{3^{\vp},4^{\vp},n\}$, it follows that $\vp(B)\not=B^{x}$, and hence $\T\not \cong \T^{'}$.
\end{example}

\section{Quotients and lifts of a triple factorisation}\label{sec-quo-lift}

\subsection{Quotients}

Let $\T=(G,A,B)$ be a triple factorisation. We study the  quotient $\T/N$ of a triple factorisation $\T=(G,A,B)$ modulo a normal subgroup $N$ of $G$ as defined in Section \ref{sec-intro}, to find conditions on $N$ under which the quotient of a nondegenerate $\T$ remains nondegenerate. We also give a sufficient condition for $\T/N$ to be nondegenerate which, in particular, allows us to assume that $G$ acts faithfully on $\Om_A$ and $\Om_B$.

For a subgroup $H\leq G$, we denote by $\ov{H}$ the corresponding subgroup $HN/N$ of $G/N$. Thus $\T/N=( \Gb, \Ab, \Bb )$ and we sometimes also denote it by $\ov{\T}$. If $\T/N$ is a nondegenerate factorisation, so is  $\T$, but the converse is not true in general as is shown by a small example in Example \ref{ex-quo}.

\begin{example}\label{ex-quo}
Let $G=\Alt(\{1,2,\ldots,5\})\times \Sym(\{6,7\})=N\times M$, say, let $B$ be the Klein $4$-subgroup of $\Alt(\{1,2,3,4\})$, and let $A=\la (1,2,3,4,5) \ra \times M$. Then $\T=(G,A,B)$ is a nondegenerate triple factorisation while  $\T/N\cong (M,M,1)$ is trivial.
\end{example}

The following lemma relates the quotient of $\T=(G,A,B)$ modulo $N$ to the lift  $\T^{'}=(G,AN,BN)$ of $\T$. It implies in particular that studying $\T^{'}$ is equivalent to studying $\T/N$.
\begin{lemma} \label{lem-quo1}
Let $\T=(G,A,B)$ be a triple factorisation and $1\not= N \lhd G$. Then  $\T/N$ is a triple factorisation. Moreover,
\begin{enumerate}
  \item[(a)] $\T/N$ is nontrivial  if and only if  $(G,AN,BN)$ is nontrivial (that is,  if and only if   both $AN$ and $BN$ are proper subgroups of $G$); and $\T/N$ is  nondegenerate  if and only if  $(G,AN,BN)$ is nondegenerate.
  \item[(b)] If  $\T$ is  nondegenerate and  $N\sse AB$, then $\T/N$ is nondegenerate.
  \item[(c)] If $N\sse \core_G(A)\core_G(B)$, then $\T/N$ is nondegenerate  if and only if  $\T$ is nondegenerate.
  \item[(d)] For  $\al=A\in \Om_A$, $\T/N$ is nondegenerate  if and only if  $\De:=\al^B \sse \Om_A$ is disjoint from at least one $N$-orbit in $\Om_A$.
\end{enumerate}
\end{lemma}
\begin{proof}
(a) This follows easily from the definitions of $\ov{A}=AN/N$ and $\ov{B}=BN/N$.\medskip

(b) Suppose that $\T$ is nondegenerate. Then $G\not =AB$. If $N\sse AB$, then $G\not =ABN$, or equivalently $\Gb\not =\Ab \ \Bb$. Therefore, $\T/N$ is nondegenerate.\medskip

(c) Suppose that $N\sse \core_G(A)\core_G(B)$. Then $N\sse AB$. Hence if $\T$  is nondegenerate, then also  $\T/N$ is nondegenerate by part (b). On the other hand, if  $\T/N$ is nondegenerate, then $(G,AN,BN)$ is nondegenerate by part (a), and hence also $\T$ is nondegenerate.\medskip

(d) Suppose that $\T/N$ is degenerate, so $\Gb=\Ab \ \Bb$, or equivalently, $G=ABN$. This implies that $BN$ is transitive on $\Om_A=[G:A]$, and hence $\Om_A=\al^{BN}=\{\de^x| \ x\in N, \ \de\in \De\}$. It follows that $\De$ contains at least one point from each $N$-orbit in $\Om_A$. Conversely, if $\De$ meets each $N$-orbit in $\Om_A$, then $\al^{BN}=\Om_A$, and hence $G=ABN$.
\end{proof}

\begin{remark} \label{rem-quo1}
\begin{enumerate}
  \item[(a)] By Lemma \ref{lem-quo1}(c), if $N=\core_G(A)\core_G(B)$, then either both or neither of $\T$ and $\T/N$ are nondegenerate. Moreover, we have $\core_{\Gb}(\Ab)=\core_{\Gb}(\Bb)=1$, so that $\Gb$ acts faithfully and transitively on both $\Om_{\Ab}:=[\Gb:\Ab]$ and $\Om_{\Bb}:=[\Gb:\Bb]$. Therefore, we may replace $\T$ by $\T/N$, and when we do so the group $G$ acts faithfully on both $\Om_A$ and $\Om_B$.
  \item[(b)] It is possible to have a nondegenerate quotient $\T/N$ without the condition $N\sse AB$ of Lemma \ref{lem-quo1}(b). For example, if $\T_i=(G_i,A_i,B_i)$ are nondegenerate, for $i=1$, $2$, then taking $G=G_1\times G_2$, $A=A_1\times A_2$, $B=B_1\times B_2$, and $N=G_2$, we have that $\T=(G,A,B)$ and $\T/N\cong \T_1$ are both nondegenerate triple factorisations, while $N\nsse AB$.
\end{enumerate}
\end{remark}

\subsection{Lifts}

Let  $\T=(G,A,B)$ be a triple factorisation, $A\leq C $ and $B\leq D$.  Recall from Section \ref{sec-intro}    that the triple factorisation $(G,C,D)$ is called a lift of $\T$. We give here an  example to illustrate that it is not possible in general to obtain a nondegenerate  lift $(G,C,D)$ of a given nondegenerate triple factorisation $\T=(G,A,B)$, with $C$, $D$ maximal in $G$.

\begin{example}\label{ex-lift}
Let $G=A_5$, $A=\la (1,2,3,4,5) \ra$, $B=\la (1,3)(2,5),(1,2)(3,5) \ra$. Then $\T=(G,A,B)$ is a nondegenerate triple factorisation.  Moreover, $A$ and $B$ are each contained in a unique maximal subgroup of $G$, namely $C\cong D_{10}$ and $D \cong A_4$, respectively, and $G=CD$. Thus the only lift  $\T^{'}=(G,C,D)$ of $\T$ with $C$, $D$ maximal is degenerate.
\end{example}

\section{Restrictions of a triple factorisation} \label{sec-res}

\begin{definition} \label{defn-res}
For a triple factorisation $\mc{T}=(G,A,B)$, and $H\leq G$, the triple  $\T|_{H}=(H,A\cap H,B\cap H)$  is called the \textit{restriction} of  $\mc{T}$ to $H$. If $\T|_{H}$ is a triple factorisation, we say that  $\T$ \textit{restricts} to $H$.
\end{definition}

We consider the following question in several contexts, in this section and again in Section \ref{sec-res-trans-normal}.

\begin{question}
Is the restriction $\T|_{H}$  a triple factorisation? If so, when is it  nontrivial, or nondegenerate?
\end{question}

\begin{remark}\label{rem-iso1}
Let  $\T:=(G,A,B)$ and $\T^{'}:=(G',A',B')$ be isomorphic triple factorisations via $\vp$, and let $H\leq G$. Then $\T$ restricts to $H$  if and only if  $\T^{'}$ restricts to $\vp(H)$. Moreover, $\T|_{H}$ is nondegenerate  if and only if  $\T^{'}|_{\vp(H)}$ is nondegenerate.
\end{remark}

\subsection{Restriction  to overgroups of $A$}\label{sec-over-A}

For any triple factorisation $\T:=(G,A,B)$ restricting to $A$ gives a trivial triple factorisation $\T|_{A}=(A,A,A\cap B)$. However,  restricting  to proper overgroups of $A$ often yields nondegenerate triple factorisations.

\begin{lemma}\label{lem-res1}
Suppose that $\T$ is a triple factorisation and  $A< H \leq G$. Then $\T|_{H}$ is a triple factorisation; moreover, $\T|_{H}$ is nontrivial  if and only if  $B\neq G$, and is nondegenerate  if and only if  $H\nsse AB$.
\end{lemma}
\begin{proof}
Let $x\in H$. Since $G=ABA$, we have $x=abc$ with $a$, $c \in A$, $b\in B$. Then $b=a^{-1}xc^{-1}\in B\cap H$, and so  $x \in A(B\cap H)A$. Thus $\T|_{H}$ is a triple factorisation. Since $H\not =A$, $\T|_{H}$ is trivial  if and only if  $B\cap H=H$,  or equivalently,  $A<H\leq B$, which is equivalent to $G=ABA=B$. If $\T|_{H}$ is degenerate, then  $H=A(B\cap H)\sse AB$. Conversely, if $H\sse AB$, then each $x\in H$ can be written  as $x=ab$ with $a\in A$, $b\in B$, and hence $b= a^{-1}x\in B\cap H$, so $H=A(B\cap H)$.
\end{proof}

The following lemma is an important special case of Lemma \ref{lem-res1}.

\begin{lemma}\label{cor-res}
Suppose that  $\T:=(G,A,B)$ is a nondegenerate triple factorisation, and that $N\unlhd G$ is such that $N\nsse A$. Then
\begin{enumerate}
  \item[(a)] $\T|_{AN}$ is a nontrivial triple factorisation, and is  nondegenerate   if and only if  $N\nsse AB$.
  \item[(b)] At least one of $\T/N$ and $\T|_{AN}$ is nondegenerate.
\end{enumerate}
\end{lemma}
\begin{proof}
Part (a) is a special case of Lemma \ref{lem-res1}. For part (b), assume that $\T/N$ is degenerate. Then by Lemma \ref{lem-quo1}(b), $N\nsse AB$, and hence $\T|_{AN}$ is nondegenerate by part (a).
\end{proof}

The proof in Lemma \ref{lem-res1} that $\T|_{H}$ is a triple factorisation  does not use the fact that $B$ is a subgroup. We give this slightly more general statement and apply it in Section \ref{sec-res-trans-normal}.

\begin{lemma}\label{lem-res4}
Suppose that $A$ is a subgroup of group $G$, and $S$ is a subset of $G$ such that $G=ASA$. Let $H$ be a subgroup of $G$. Then $A\leq H$  if and only if  $H=A(S\cap H)A$.
\end{lemma}
\begin{proof}
Clearly $H=A(S\cap H)A$ implies $A\leq H $. Conversely, suppose that $A\leq H$. Repeating the proof of Lemma \ref{lem-res1}, we get $H= A(S\cap H)A$.
\end{proof}

\subsection{Restriction to  normal subgroups}

If $\T:=(G,A,B)$ and $\T^{'}:=(G',A',B')$ are triple factorisations with $G=G'$, $A=A'$, $B\cong B'$, and if $N$ has index $2$ in $G$, it is possible for $\T$ but not $\T^{'}$ to restrict to $N$.

\begin{lemma}\label{lem-res3}
There exist infinitely many triple factorisations $\T=(G,A,B)$ and  $\T^{'}:=(G,A,B^{x})$, where $x\in G$, such that $G$ has a subgroup $N$ of index $2$ and $\T$ restricts to $N$ but $\T^{'}$ does not restrict to $N$.
\end{lemma}

\begin{proof}
Consider the factorisations $\T=(G,A,B)$ and $\T^{'}=(G,A,B')$ in Example~\ref{ex-part(b)}, and let $N=\Alt(\Om)$. We claim that $\T$ restricts to $N$ but $\T^{'}$ does not.

Now $\T|_{N}=(N,G_n\cap N, \la (3,4,n) \ra)$ while $\T^{'}|_{N}=(N, G_n\cap N, B')$ where  $B'= \la (3,4,1)\ra$. In the latter case both $A\cap N$ and $B'$ fix the point $n$ and hence $\T^{'}$ does not restrict to $N$. On the other hand, if $g\in N$ fixes $n$, then $g \in A\cap N$, while if $j=n^g\not=n$, then there exists $c\in A\cap N$ such that $3^c=j$ and hence setting $b=(3,4,n) \in B\cap N$, the element $a=gc^{-1}b^{-1}$ fixes $n$ and so lies in $A\cap N$. Thus $g=abc \in (A\cap N)(B\cap N)(A\cap N)$, and it follows that $\T$ restricts to $N$.
\end{proof}

We explore this problem further in Section
\ref{sec-res-trans-normal}.
\section{Wreath products of triple factorisations} \label{sec-wreath}

\noindent  In this section, we introduce  a wreath product construction for triple factorisations, and study its properties. In Subsection \ref{sec-reduction-thm}, we prove  Theorem \ref{thm:main}.

Let $\T_0=(G_0,A_0,B_0)$ and $\T_1=(G_1,A_1,B_1)$ be triples with $A_i$, $B_i$ subgroups of $G_i$, for $i=0$, $1$ (not necessarily triple factorisations), and let $\De$ and $\Si$ denote the sets of right cosets of $A_0$ and $A_1$ in $G_0$ and $G_1$, respectively. Recall from Section \ref{sec-embed} that the set $\Fun(\Si,G_0)$ consisting of all functions from $\Si$ to $G_0$ forms a group under pointwise multiplication, that is to say, $(fg)(\De)=f(\De)g(\De)$,  and the wreath product $G_0\wr G_1$ is the semidirect product $\Fun(\Si,G_0)\semi G_1$ acting on $\De\times \Si$ via
\[(\de,\De)^{f\si}=(\de^{f(\De)},\De^{\si}) \]
for all $\de \in \De$, $\De\in \Si$, and $f\si\in G_0\wr G_1$. Note that $G_0\wr G_1$ can also be identified with  $G_0^{\ell}\semi G_1$, where $|\Si|=\ell$, and for all $\si \in G_1$ and $(f_1,...,f_\ell)\in G_0^{\ell}$, we have
\[\si(f_1,...,f_\ell)\si^{-1}=(f_{1^\si},...,f_{\ell^\si}).\]
We define a wreath  product of $\T_0$ and $\T_1$ as follows.

\begin{definition} \label{defn-wreath}
Let $\T_0=(G_0,A_0,B_0)$ and $\T_1=(G_1,A_1,B_1)$ be triples as above, and let $\De=[G_0:A_0]$, $\Si=[G_1:A_1]$ and $\ell=|\Si|$. Then the \emph{wreath product} of $\T_0$ and $\T_1$ is defined as  $\T_0\wr \T_1=(G_0\wr G_1, \h{A},\h{B})$, where $\h{A}=(A_0\tm G_0^{\ell-1})\semi A_1$ and $\h{B}=B_0\wr B_1$.
\end{definition}

\begin{lemma} \label{lem-wr1}
Let $\h{A}$ and $\h{B}$ be  as in Definition \ref{defn-wreath}, and let $W=G_0\wr G_1$. Then
$\h{A}\h{B}= (A_0B_0\tm G_0^{\ell-1}) (A_1B_1)$, $\core_{W}(\h{A})=\core_{G_0}(A_0)\wr \core_{G_1}(A_1)$ and $\core_{W}(\h{B})=\core_{G_0}(B_0)\wr \core_{G_1}(B_1)$.
\end{lemma}
\begin{proof}
Let $f\si:=(a_1, \ldots,a_\ell)\si \in \h{A}$ and $g\tau:=(b_1, \ldots,b_\ell)\tau \in \h{B}$, where $a_1\in A_0$, $a_i\in G_0$ for $i\geq 2$, $b_i\in B_0$ for $i\geq 1$, $\si\in A_1$ and  $\tau \in B_1$. Since $\si \in A_1$ and  $b_{1^{\si}}\in B_0$, we have
\begin{align*}
    (f\si)(g\tau)&=[(a_1, \ldots,a_\ell)\si][(b_1, \ldots,b_\ell)\tau]\\
                 &= [(a_1,\ldots,a_\ell)(b_{1^{\si}},\ldots,b_{\ell^{\si}})] (\si \tau) \\
                 &=(a_1b_{1^{\si}},a_2b_{2^{\si}},\ldots,a_\ell b_{\ell^{\si}}) (\si \tau)\\
                 &\in (A_0B_0\tm G_0^{\ell-1}) (A_1B_1).
\end{align*}
Therefore $\h{A}\h{B}\sse (A_0B_0\tm G_0^{\ell-1}) (A_1B_1)$. The converse is similar and easier. Finally the assertions about $\core_{W}(\h{A})$ and $\core_{W}(\h{B})$ follow since $G_1$ is transitive on $\Si$.
\end{proof}

\begin{lemma}\label{lem-wr2}
Let $\T_0$ and $\T_1$ be as in Definition \ref{defn-wreath}. If $\T_0$ and $\T_1$ are both triple factorisations, then $\T_0\wr \T_1$ is a triple factorisation. Moreover, $\T_0\wr \T_1$ is trivial (nondegenerate)  if and only if  at least one of the $\T_0$ and $\T_1$ is trivial (nondegenerate, respectively).
\end{lemma}
\begin{proof}
Let  $\h{A}:=(A_0\tm G_0^{\ell-1})\semi A_1 $, $\h{B}:=B_0\wr B_1$, and $W:=G_0\wr G_1$, so $W=G_0^{\ell}\semi G_1$. Suppose that $f\si \in W$, where  $f:=(f_1,\ldots, f_\ell)\in G_0^\ell$ and $\si \in G_{1}$. Since $\T_0$ is a triple factorisation, there exist $a_i, c_i \in A_0$, $b_i \in B_0$ such that $f_i=a_ib_ic_i$, for $i=1,\ldots, \ell$. Also since $\T_1$ is a triple factorisation, there exist $\de$, $\nu \in A_1$ and $\tau \in B_1$ such that $\si=\de \tau \nu$. Let $\h{a}:=(a_1,\ldots, a_\ell) \de \in \h{A}$, $\h{b}:=(b_{1^{\de^{-1}}},\ldots, b_{\ell^{\de^{-1}}} )\tau \in \h{B}$, $\h{c}:=(c_{1^{(\de\tau)^{-1}}},\ldots,$ $c_{\ell^{(\de\tau)^{-1}}} ) \nu \in \h{A}$. Then
\begin{align*}
    f\si&=(a_1b_1c_1, ...,a_\ell b_\ell c_\ell ) \de\tau\nu\\
     &= (a_1,\ldots,a_\ell)(b_1,...,b_\ell )\de\tau[(c_{1^{(\de\tau)^{-1}}},...,c_{\ell^{(\de\tau)^{-1}}} ) \nu]\\
     &= [(a_1,\ldots,a_\ell)\de][(b_{1^{\de^{-1}}},...,b_{\ell^{\de^{-1}}} )\tau] [(c_{1^{(\de\tau)^{-1}}},...,c_{\ell^{(\de\tau)^{-1}}} ) \nu]\\
     &=\h{a}\h{b}\h{c} \in \h{A}\h{B}\h{A}.
\end{align*}
Therefore, $W=\h{A}\h{B}\h{A}$, so $\T_0\wr \T_1$ is a triple factorisation.

By the definition, $\T_0\wr T_1$ is trivial  if and only if  $W=\h{A}$ or $W= \h{B}$. Since $\h{A}:=(A_0\tm G_0^{\ell-1})\semi A_1$ and $\h{B}:=B_0\wr B_1$, $\T_0\wr \T_1$ is trivial  if and only if  either (a)  $G_0=A_0$ and $G_1=A_1$, or (b) $G_0=B_0$ and $G_1=B_1$. Each of the conditions (a) or (b) holds  if and only if  at least one of  $\T_0$ and $\T_1$ is trivial. By Lemma \ref{lem-wr1}, we have  $\h{A}\h{B}=(A_0B_0\tm G_0^{\ell-1})A_1B_1$, and hence $\T_0\wr \T_1$ is nondegenerate  if and only if  $G_0\not=A_0B_0$ or $G_1\not=A_1B_1$, or equivalently, at least one of $\T_0$ and $\T_1$ is nondegenerate.
\end{proof}

\subsection{Proof of the reduction Theorem \ref{thm:main}} \label{sec-reduction-thm}

Suppose that $\T:=(G,A,B)$ is a nondegenerate triple factorisation of $G$, and that $H$ satisfies $A<H<G$ with $\core_G(A)=1$. Let $\Om:=\Om_{A}=[G:A]$ with $G$ acting by right multiplication. Set $\al:=A \in \Om$,  $\De:=\al^H$, and let $\Si:=\{\De_1, \ldots, \De_\ell\}$ be the block system determined by $H$, where $\De=\De_1$ and $|G:H|=\ell$. Let now $G_0$ and $G_1$ be the permutation groups induced by $H$ on $\De$ and $G$ on  $\Si$, respectively. Let $A_1=H^{\Si}$, $B_1=B^{\Si}$, $A_0=A^{\De}$ and $B_0=(B\cap H)^{\De}$. By the Embedding Theorem \ref{thm:embed}, we may assume that $\Om=\De\times \Si$ and that $G$ is a subgroup of $W=G_0\wr G_1$ such that $A=\h{A}\cap G$ where  $\h{A}=(A_0\times G_0^{\ell-1}) \semi A_1$. Moreover, if $B_1$ is transitive, then by Theorem \ref{thm:embed} we may also assume that $B\leq \h{B}$ where $\h{B}=B_0 \wr B_1$. Define
\begin{equation} \label{defn-wr1}
\T_0:=(G_0,A_0,B_0),\quad \T_1:=(G_1,A_1,B_1), \quad \T_0\wr \T_1:=(W,\h{A},\h{B})
\end{equation}
with $\T_0\wr T_1$ as in Definition \ref{defn-wreath}.

\begin{lemma}\label{lem-wr4}
Let $\h{A}$ and $\h{B}$ be as above, and $\h{H}:=G_0\wr A_1$. Then $\core_{\h{H}}(\h{A})=(1\tm G_0^{\ell-1})\semi A_1$ and $\core_{W}(\h{A})=1$.
\end{lemma}
\begin{proof}
By the definition of $G_0$ and $G_1$, it follows that $\core_{G_0}(A_0)=1$ and $\core_{G_1}(A_1)=1$, and hence by Lemma \ref{lem-wr1}, $\core_{W}(\h{A})=1$. Let $\h{a}:=(a_1, \ldots,a_\ell)\si \in \h{A}$ and $\h{h}:=h\nu:=(h_1, \ldots,h_\ell)\nu \in \h{H}$, where $a_1\in A_0$, $a_i\in G_0$ for $i\geq 2$, $h_i\in G_0$ for $i\geq 1$, $\si$, $\nu \in A_1$. Since  $(h\nu)^{-1}={(h^{\nu})}^{-1}\nu^{-1}$, $1^\si=1$ and  $1^\nu=1$, we have
\begin{align*}
\h{h}^{-1}\h{a}\h{h}&=((h^{\nu})^{-1}\nu^{-1})\h{a}(h\nu)\\
                &=[(h_{1^{\nu^{-1}}}^{-1},\ldots,h_{\ell^{\nu^{-1}}}^{-1})\nu^{-1}][(a_1, \ldots,a_\ell)\si][(h_1, \ldots,h_\ell)\nu]\\
                &= [(h_{1}^{-1},h_{2^{\nu^{-1}}}^{-1},\ldots,h_{\ell^{\nu^{-1}}}^{-1})(a_1, a_{2^{\nu^{-1}}} \ldots,a_{\ell^{\nu^{-1}}})\\
                &\hspace{3cm} (h_1, h_{2^{\nu^{-1}\si}},\ldots ,h_{\ell^{\nu^{-1}\si}})][\nu^{-1} \si \nu] \\
                &= (h_{1}^{-1}a_1h_1,h_{2^{\nu^{-1}}}^{-1}a_{2^{\nu^{-1}}}h_{2^{\nu^{-1}\si}},\ldots ,h_{\ell^{\nu^{-1}}}^{-1}a_{\ell^{\nu^{-1}}}h_{\ell^{\nu^{-1}\si}})[\nu^{-1} \si \nu].
\end{align*}
Therefore, $(\h{A})^{\h{h}}=(A_0^{h_1}\tm G_0^{\ell-1})\semi A_1$, and hence $\core_{\h{H}}(\h{A})=(\core_{G_0}(A_0)\tm G_0^{\ell-1})\semi A_1=(1 \tm G_0^{\ell-1})\semi A_1$.
\end{proof}

\begin{proposition}\label{prop-wr1}
Suppose that $\T:=(G,A,B)$ is a triple factorisation for $G$ with $\core_G(A)=1$, and suppose that $A<H<G$. Let $\T_0$, $\T_1$ and $\T_0\wr \T_1$ be as in \emph{(\ref{defn-wr1})}. Then
\begin{enumerate}
  \item[(a)] $\T_0$ is a triple factorisation, and  is nontrivial  if and only if  $H\not = \core_H(A)(B\cap H)$, and nondegenerate  if and only if  $H\nsse AB$.
  \item[(b)]  $\T_1$ is  a triple factorisation, and is the quotient of the lift $\T^{'}=(G,H,B)$ of $\T$ modulo  $\core_G(H)$. Moreover,  $\T_1$ is nontrivial  if and only if  $G\not =\core_G(H)B$, and nondegenerate  if and only if  $G\not =HB$.
  \item[(c)] $\T_0\wr \T_1$ is a triple factorisation, and is nontrivial  if and only if  either $H\not = \core_H(A)(B\cap H)$, or $G\not =\core_G(H)B$. Also $\T_0\wr \T_1$ is nondegenerate  if and only if  either $H\nsse AB$, or $G\not=HB$.
  \item[(d)] If $\T$ is nondegenerate, then also $\T_0\wr \T_1$ is nondegenerate.
\end{enumerate}
\end{proposition}
\begin{proof}
\noindent (a) Since $A<H$, by Lemma \ref{lem-res1}, $\T|_{H}=(H,A,B\cap H)$ is  a nontrivial triple factorisation, and is nondegenerate  if and only if  $H\nsse AB$. Let $N=\core_{H}(A)$, so that $\T_0=(\T|_{H})/N$. It follows from Lemma \ref{lem-quo1}(a) that $\T_0$ is a triple factorisation for $G_0$, and $\T_0$ is nontrivial  if and only if   the lift $(H,AN,(B\cap H)N)$ of $\T|_{H}$  is nontrivial. Since $N\subset A$, thus latter holds  if and only if  $(B\cap H)N\not= H$. Moreover, by Observation \ref{modcore},  $\T_0$ is nondegenerate  if and only if  $\T|_{H}$ is nondegenerate, or equivalently,  if and only if  $H\nsse AB$.\medskip

(b) Now let $N=\core_{G}(H)$ and note that $H(BN)=HB$. Since $A<H$, the triple  $\T^{'}:=(G,H,B)$ is a lift of  $\T$, and by the definition of a quotient triple factorisation,  $\T_1=\T^{'}/N$ (and is also a lift of $\T/N=(G/N,AN/N,BN/N)$). By Lemma \ref{lem-quo1}(a), $\T_1$ is a triple factorisation, and is nontrivial  if and only if  $G\not=BN$. Moreover by Observation \ref{modcore}, $\T_1$ is nondegenerate  if and only if  $\T^{'}$ is nondegenerate, that is to say,  if and only if  $G\not=H(BN)$. As we observed above this is equivalent to $G\not=HB$.\medskip

(c) By parts (a) and (b), $\T_0$ and $\T_1$ are triple factorisations, and so by Lemma \ref{lem-wr2}, $\T_0\wr \T_1$ is a triple factorisation. Also, by Lemma \ref{lem-wr2}, $\T_0\wr T_1$ is nontrivial  if and only if  at least one of the $\T_0$ and $\T_1$ is nontrivial and hence, by parts (a) and (b),  if and only if  either $H\not = \core_H(A)(B\cap H)$, or $G\not =\core_G(H)B$. By Lemma \ref{lem-wr2},  $\T_0\wr \T_1$ is nondegenerate  if and only if  at least one of $\T_0$ and $\T_1$ is nondegenerate and hence, by parts (a) and (b),  if and only if  either $H\nsse AB$, or $G\not=HB$.\medskip

(d) Let $\T$ be nondegenerate. By part (c), $\T_0\wr \T_1$ is nondegenerate  if and only if  either  $H\nsse AB$, or $G\not=HB$. Suppose that $G=HB$. If also $H\sse AB$, then $G=HB\sse ABB=AB$, and so $G=AB$, or equivalently, $\T$ is degenerate, which is a contradiction. Therefore in this case  $H\nsse AB$. Thus at least one of $G\neq HB$ or $H\nsse AB$ holds, and hence by  part (c),  $\T_0\wr \T_1$ is nondegenerate.
\end{proof}

We define several subgroups of $W=G_0\wr G_1$ that correspond to the subgroups $H$, $\core_H(A)$ and $\core_G(H)$ of $G$ that occur in Proposition \ref{prop-wr1}. Recall that  $\h{A}=(A_0\times G_0^{\ell-1}) \semi A_1$ and $\h{B}=B_0 \wr B_1$.
\begin{equation}\label{Wsubgroups}
\h{H}:=G_0\wr A_1, \quad \h{K}:=(1 \tm G_0^{\ell-1})\semi A_1, \quad \h{N}:=G_0^\ell.
\end{equation}

\begin{lemma} \label{lem-wr3}
Suppose that $\T=(G,A,B)$ is a nondegenerate triple factorisation  with $\core_G(A)=1$, and $A<H<G$. Let $\T_0$, $\T_1$, $\T_0\wr \T_1$ and $\ell$ be as in \emph{(\ref{defn-wr1})}, and assume that $G\leq W=G_0\wr G_1$ so that $A=\h{A}\cap G$. Also  if $B_1$ is transitive, assume that $B\leq \h{B} \cap G$. Let $\h{N}$, $\h{H}$ $\h{K}$ be as in \emph{(\ref{Wsubgroups})}. Then the following hold.
\begin{enumerate}
  \item[(a)] $\h{K}=\core_{\h{H}}(\h{A})$, $\h{H}/\h{K}\cong G_0$, and  $\T^{'}/\h{K}\cong \T_0$, where $\T^{'}=(\T_0\wr \T_1)|_{\h{H}}$.
  \item[(b)] $\h{N}=\core_{W}(\h{H})$ and $(\T_0\wr \T_1)/\h{N}\cong \T_1$.
  \item[(c)] If $\T_1$ is degenerate, then $\T_0$ and $\T_0\wr \T_1$ are both nondegenerate and  $(\T_0\wr \T_1)|_{G}=(G,A,\h{B}\cap G)$ is a nondegenerate lift of $\T$.
\end{enumerate}
\end{lemma}
\begin{proof}
(a) By (\ref{Wsubgroups}), we have  $\h{A}< \h{H}$, and by Lemma \ref{lem-res1},  $(\T_0\wr\T_1)|_{\h{H}}=(\h{H},\h{A},\h{B}\cap\h{H})$ is a nontrivial triple factorisation. Note that $\h{B}\cap\h{H}=B_0\wr (A_1\cap B_1)$, and also by Lemma \ref{lem-wr4} that  $\h{K}=\core_{\h{H}}(\h{A})$. Consider  the projection map $\vp:\h{H}\longrightarrow G_0$ given by $(h_1,\ldots,h_\ell)\si\mapsto h_1$ (for $h_i\in G_0$ and $\si\in A_1$). Then $\vp$ is a epimorphism with kernel $\h{K}$, $\vp(\h{A})= A_0$, $\vp(\h{B})= B_0$, and $\vp(\h{H})\cong \h{H}/\h{K}$. Hence $\vp$ induces an isomorphism from the quotient of $(\T_0\wr \T_1)|_{\h{H}}$ modulo $\h{K}$ to $\T_0$.\medskip

(b) By (\ref{Wsubgroups}), $(G_0\wr G_1)/\h{N}\cong G_1$, $\h{A}\h{N}=G_0\wr A_1=\h{H}$ and $\h{B}\h{N}=(B_0\wr B_1)G_0^{\ell}=G_0\wr B_1$. Thus $W/\h{N}\cong G_1$, $\h{A}\h{N}/\h{N}\cong A_1$ and $\h{B}\h{N}/\h{N}\cong B_1$, and hence the natural projection map $W\longrightarrow W/\h{N}$ defines an isomorphism from $(\T_0\wr \T_1)/\h{N}$ to $\T_1$.\medskip

(c) Suppose that $\T_1$ is degenerate, so by Proposition \ref{prop-wr1}(b), $G=HB$, or equivalently,  $B_1$ is transitive. Thus by our assumption $B\leq \h{B}\cap G$. Hence the restriction of $\T_0\wr \T_1=(W,\h{A},\h{B})$ to $G$ is $\T^{'}:=(G,A, \h{B}\cap G)$, a lift of $\T$. The subset $A(\h{A}\cap G)\sse \h{A}\h{B}$, and so by Lemma \ref{lem-wr1},   $A(\h{A}\cap G)\sse (A_0B_0\tm G_0^{\ell-1})\semi A_1B_1$. We will construct an element $h$ lying in $G$ but not in $(A_0B_0\times G_0^{\ell-1})\semi A_1B_1$, proving that $ A(\h{A}\cap G)\neq G$, and hence that $\T^{'}$ is nondegenerate. By Proposition \ref{prop-wr1}(d) and Lemma \ref{lem-wr2}, and since $\T_1$ is degenerate, it follows that $\T_0$ and $\T_0\wr \T_1$ are both nondegenerate. Thus $G_0\not=A_0B_0$. Let $g\in G_0\sm A_0B_0$. Since $G_0=H^{\De}$, $g=h^{\De}$ for some $h\in H=(G_0^{\ell}\semi A_1)\cap G$. This element $h$ is therefore equal to $(g,h_2,\ldots, h_{\ell})\si$ for some $h_2$,...,$h_{\ell} \in G_0$ and $\si \in A_1$, and therefore $h$ lies in $G$ but not in $(A_0B_0\times G_0^{\ell-1})\semi A_1B_1$.
\end{proof}

\noindent \emph{Proof of Theorem \ref{thm:main}. } Suppose that $\T=(G,A,B)$ is a nondegenerate triple factorisation with $\core_G(A)=1$, that $A<H<G$, and that $G\leq G_0\wr G_1$ with $G_0$ and $G_1$ as above. Let $\T_0$, $\T_1$ and $\T_0\wr \T_1$ be as in (\ref{defn-wr1}). By the definitions of $\T_0$ and $\T_1$, we have that  $\T_0$ is a quotient of $\T|_{H}$, and $\T_1$ is the quotient modulo $\core_G(H)$ of the lift $(G,H,B)$ of $\T$. By Proposition \ref{prop-wr1}, all three are triple factorisations, and $\T_0\wr \T_1$ is nondegenerate. Finally if $\T_1$ is degenerate, then  $B_1$ is transitive and by Theorem \ref{thm:embed} we may assume that $B\leq \h{B}$. Then the assertions of Theorem \ref{thm:main}(b) all follow from Lemma \ref{lem-wr3}. \hfill $\Box$

\begin{corollary} \label{cor-wr4}
If $B$ is a  maximal subgroup of $G$ and $B_1=B^{\Si}$ is transitive, then  $(\T_0\wr \T_1)|_{G}=\T$.
\end{corollary}
\begin{proof}
Since $B_1$ is transitive, by Theorem \ref{thm:main}(b), $B\leq \h{B}\cap G$ and $(\T_0\wr \T_1)|_{G}=(G,A,D)$ is nondegenerate, for some $D$ such that $B\leq D\leq G$. In particular, $D\neq G$. Then since $B$ is maximal, it follows that $D=B$, so $(\T_0\wr \T_1)|_{G}=\T$.
\end{proof}

\noindent \emph{Proof of Corollary \ref{cor:main}. }

Let $\T=(G,A,B)$ be a nondegenerate triple factorisation with  $\core_G(A)=1$.

(a) Suppose first that $A$ is maximal in $G$. Then $G$ acts faithfully and primitively on $\Om_A$. If $X=\core_G(B)\neq 1$, then since $G$ is primitive, $X$ is transitive on $\Om_A$. Thus $G=AX\sse AB$, contradicting the nondegeneracy of $\T$. Hence $X=1$.\medskip

(b) Now let $A=H_1<H_2<...<H_r=G$ with $H_i$ maximal in $H_{i+1}$ for $1\leq i<r$. Note that $G=H_r=BH_r$. Let $j$ be minimal such that $G=BH_j$. Since $G\neq BA$, we have $2\leq j\leq r$. Let $H=H_j$, $K=H_{j-1}$. By the minimality of $j$, $G\neq BK$. However $G=ABA=KBK$, and hence $\T^{'}=(G,K,B)$ is a nondegenerate lift of $\T$. By Theorem \ref{thm:main} applied to $\T^{'}$ modulo $N=\core_G(K)$, we find that $\T_0=(H,K,BN\cap H)$ modulo $M=\core_H(K)$ is nondegenerate with $K/M$ maximal in $H/M$ and $\core_{H/M}(K/M)=1$. By part (a) above, $\T_0/M$ is primitive. \hfill $\Box$

\section{Restriction to transitive normal subgroups} \label{sec-res-trans-normal}

Let $\T=(G,A,B)$ be a nondegenerate triple factorisation. By Observation \ref{modcore} and Theorem \ref{thm:main}, it is important to study such triple factorisations in which $A$ is maximal and core free in $G$, so that $G$ is faithful and primitive on $\Om_A=[G:A]$. In this case also $B$ is core-free in $G$ by Corollary \ref{cor:main}. The O'Nan-Scott Theorem (see \cite[Chapter IV]{DixMort96}) describes various types of finite primitive permutation groups $G$ identifying the types by the structure and permutation action of their socles (the \emph{socle} $\Soc (G)$ of a group $G$ is the product of the minimal normal subgroups of $G$). Thus it is natural to seek conditions under which a primitive triple factorisation $\T=(G,A,B)$ restricts to $\Soc(G)$. The examples given for Lemma \ref{lem-res3} suggest that this may be a difficult problem. In Subsection \ref{sec-resconditions}, we give a rather technical sufficient condition for restriction. Before that we prove a simple result which explores the role of normal subgroups for primitive $\T$.
\begin{lemma} \label{lem-prim1}
Suppose that $\T=(G,A,B)$ is a triple factorisation.

\begin{enumerate}
 \item[(a)] If $1\not = N\lhd G$ and $N$ is transitive on $\Om_A$, then $\T$ is nondegenerate  if and only if  $N\nsse AB$.
 \item[(b)] If $\core_G(A)=1$ and $A$ is maximal in $G$, then $\T$ is nondegenerate  if and only if  the set $AB$ contains no nontrivial normal subgroup of $G$.
\end{enumerate}
\end{lemma}
\begin{proof}
Since $N$ is transitive on $\Om_A$, $G=AN$. Thus $N\sse AB$  if and only if  $G=AN\sse AB$, and the latter holds  if and only if  $\T$ is degenerate. This proves (a) and part (b) follows immediately.
\end{proof}

Thus a triple factorisation $\T=(G,A,B)$ is primitive  if and only if  $A$ is maximal in $G$ and $AB$ contains no nontrivial normal subgroup of $G$.

\subsection{Sufficient conditions for restriction to a transitive normal subgroup} \label{sec-resconditions}

Let $G$ be a finite group with subgroups $A$, $B$ and let $N$ be a normal subgroup of $G$ such that $G=AN$.

Let $A_0=A\cap N$.  Then  $A=\la A_0, \si_1, \ldots, \si_r \ra$ for some elements $\si_i \in A\sm A_0$.  Hence $G=\la N, \si_1, \ldots, \si_r \ra$. Let $A_\si:=\la \si_1, \ldots, \si_r \ra$. Then $G=AN=A_{\si}N$ and $A=A_0A_{\si}=A_{\si}A_0$ since $A_0\unlhd A$. Moreover, if $B_0=B\cap N$, then $B=\la B_0, \tau_1, \ldots, \tau_s \ra$ for some elements $\tau_j\in B\sm B_0$. Since $G=A_{\si}N$, each $\tau_j= \de_j n_j$ for some $\de_j \in A_{\si}$, $n_j\in N$. Now set $B_\tau:=\la \tau_1, \ldots, \tau_s \ra$. Then $B=B_0B_\tau$ because $B_0\unlhd B$. We will define a subset $B_{\si}$ of $G$ such that
\begin{equation} \label{eq3-3}
G=ABA \  \Longleftrightarrow \ G=A_0(B_{\si}A_{\si})A_0.
\end{equation}
Note that $B_\si A_\si$ will not in general be a subgroup.

Suppose that $G=ABA$. If $x\in G$, then $x=abc$, where $a,c \in A$ and $b\in B$. There exist $a_0, c_0 \in A_0$ and $\lambda, \lambda'\in A_\si$ such that $a=a_0\lambda$ and $c=\lambda' c_0$. Then
\begin{align*}
    x&=abc\\
     &=(a_0\lambda )b(\lambda'c_0) \\
     &=a_0 (\lambda b \lambda^{-1}) (\lambda \lambda') c_0\\
     &=a_0 b^{\lambda^{-1}}  \lambda'' c_0,\\
\end{align*}
where $\lambda'':=\lambda \lambda'\in A_\si$. On the other hand, if $x\in A_0B^{\lambda^{-1}} A_{\si}A_0$ for some $\lambda \in A_0$, then $x=a_0b^{\lambda^{-1}}\lambda''c_0$, for some $a_0,c_0\in A_0$, $b\in B$, $\lambda''\in A_{\si}$ and it follows that $x=abc$, where $a=a_0\lambda \in A$, $c=\lambda^{-1}\lambda'' c_0 \in A$, that is, $x\in ABA$. Therefore,

\begin{equation}\label{eq3-2}
G=ABA \  \Longleftrightarrow \    G=\bigcup_{\lambda\in A_\si}A_0B^{\lambda} A_{\si} A_0.
\end{equation}
Define $B_{\si}=\cup_{\lambda\in A_\si}B^{\lambda} $. Then (\ref{eq3-3}) holds.\\
Moreover, since $B=B_0B_{\tau}$, we have $B_{\si}=\cup_{\lambda \in A_{\si}}B_0^{\lambda} B_{\tau}^\lambda$, and also
\begin{equation}\label{eq3-4}
G=ABA \  \Longleftrightarrow \  G=A_0(\cup_{\lambda \in A_{\si}}B_0^{\lambda} B_{\tau}^\lambda )A_{\si} A_0.
\end{equation}

\begin{lemma}\label{lem-res-suf}
Let $\mc{T}=(G,A,B)$ be a triple factorisation for $G$, and let $N$ be a nontrivial normal subgroup of $G$ such that $G=AN$. Suppose furthermore that $A_0$,  $A_\si$, $B_0$, $B_\tau$ and $B_\si$ are as above. Then
\begin{enumerate}
  \item[(a)] $N=A_0(B_{\si}A_{\si}\cap N) A_0$.
  \item[(b)] If $B_{\si}A_{\si}\cap N\sse B_0A_0$, then $\T$ restricts to $N$.
  \item[(c)] If $A_\si \sse N_{G}(B_0)$ and $B_\tau \sse B_0A$, then $B_\si A_\si \cap N\sse B_0A_0$, so $\T$ restricts to $N$.
\end{enumerate}
\end{lemma}
\begin{proof}
(a) By (\ref{eq3-4}), we have $G=A_0(B_\si A_\si)A_0$. Since $A_0=A\cap N\leq N$, applying Lemma \ref{lem-res4} with $(G,A,S,H)=(G,A_0,B_\si A_\si,N)$, we have $N=A_0(B_\si A_\si \cap N)A_0$ and part (a) follows. \medskip

(b) If $B_{\si}A_{\si}\cap N\sse B_0A_0$, then by part (a), $N=A_0(B_{\si}A_{\si}\cap N) A_0\sse A_0(B_0A_0)A_0=A_0B_0A_0\sse N$. Thus $N=A_0B_0A_0$, or equivalently, $\T$ restricts to $N$. \medskip

(c) Suppose that $A_\si \leq N_G(B_0)$ and $B_\tau \sse B_0A_0$. Then $B_\si=\cup_{\lambda \in A_\si}(B_0B_\tau)^\lambda=B_0(\cup_{\lambda \in A_\si}B_\tau^\lambda)$, so that
\[B_\si A_\si=B_0(\cup_{\lambda \in A_\si}(B_\tau A_\si)^\lambda).\]
Since $B_\tau \sse B_0A$, we have, for each $\lambda \in A_\si$, that  $(B_\tau A_\si)^\lambda\sse (B_0A)^\lambda=B_0A$. So $B_\si A_\si \sse B_0A$, and hence  $(B_\si A_\si)\cap N\sse (B_0A) \cap N= B_0(A \cap N)= B_0A_0$. Therefore by part (b), $\T$ restricts to $N$.
\end{proof}

Using some of the examples from Example \ref{ex-part(b)}, we illustrate how Lemma \ref{lem-res-suf} can be used to prove that triple factorisations restrict to certain normal subgroups.

\begin{example}
Suppose  that $G=\Sym(\Om)$ acting on  $\Om=\{1,\ldots, n\}$ with $n \geq 5$, $N=\Alt(\Om)$, and $\al=n \in\Om$. Let $A=G_\al\cong S_{n-1}$, so that $A_0:=A \cap N=N_\al \cong A_{n-1}$, and $A=\la A_0, \si \ra$, where $\si$ is any transposition fixing $\al$. As in Example \ref{ex-part(b)}, let $B=\la (1,2),(3,4,n)\ra$, so $B_0:=B\cap N=\la (3,4,n) \ra$. To apply Lemma \ref{lem-res-suf}, we choose $\si\in N_G(B_0)$, for example $\si=(1,2)$, so that $A_\si=\la (1,2)\ra$. Also if we choose $\tau =(1,2)=\si$, then $B=\la B_0,\tau \ra$, $B_\si=B$,  and $B_{\tau}=\la \tau \ra \sse B_0A_{\si}$. Thus by Lemma \ref{lem-res-suf}(c), $\T$ restricts to $N$. Note also that $A_\si=\la \si \ra\cong Z_2$ and $B_\si=B\cup B^{\si}=B$.
\end{example}


\noindent \textbf{Acknowledgment:} Both authors are grateful to John Bamberg for helpful advice. The first author is grateful for support of an International Postgraduate Research Scholarship (IPRS).
This project forms part of Australian Research Council Discovery Grant DP0770915. The second author is supported by Australian Research Federation Fellowship FF0776186.

\bibliographystyle{amsplain}
\def\cprime{$'$}
\providecommand{\bysame}{\leavevmode\hbox to3em{\hrulefill}\thinspace}
\providecommand{\MR}{\relax\ifhmode\unskip\space\fi MR }
\providecommand{\MRhref}[2]{%
  \href{http://www.ams.org/mathscinet-getitem?mr=#1}{#2}
}
\providecommand{\href}[2]{#2}

\end{document}